\documentclass[11pt]{amsart} 
\pdfmapfile{=newtx.map}
\pdfmapfile{=cm-super-t1.map}
\usepackage[parfill]{parskip}
\usepackage[T1]{fontenc}
\usepackage{amsthm}
\usepackage{textcomp}
\usepackage{blkarray}
\usepackage[scaled=.85]{beramono}% a typewriter font must be defined
\usepackage[leqno]{amsmath} %use eq no on left
\usepackage[cmintegrals,cmbraces,libertine]{newtxmath}
\usepackage{bm}
\usepackage[bb=boondox]{mathalfa}
\usepackage{bm}% load after all math to give access to bold math

\arraycolsep=1.4pt\def\arraystretch{1.4}

\newcommand{\bP}{{\varmathbb{P}}}

\newcommand{\bK}{{\varmathbb{K}}}
\def\Mid{\,\vert\,}

\usepackage{graphicx}
\usepackage{calc}
\usepackage{epstopdf}

\usepackage[paperwidth=7in, paperheight=10in, textwidth=4.75in, textheight=8in, includehead=true, bindingoffset=.5in]{geometry}

\newcommand{\beq}{\begin{equation}}
\newcommand{\eeq}{\end{equation}}
\numberwithin{equation}{section}

\newtheorem{theorem}[equation]{Theorem}

\newtheorem{lemma}[equation]{Lemma}

\theoremstyle{remark}
\newtheorem{srem}[equation]{}
\newtheorem{sremark}[equation]{}

\newtheorem{remark}[equation]{Remark}

\theoremstyle{definition}

\newtheorem{definition}[equation]{Definition}

\title{\textsf{Schur Inequalities}}%\\[15pt]{\em proof techniques and intuition}}
\author{S. Gill Williamson}
\thanks{Department of Computer Science and Engineering, 
University of California San Diego; \url{http://cse.ucsd.edu/~gill}.
}
\date{}                                           % Activate to display a given date or no date
\usepackage{setspace}
\usepackage[pdftex,colorlinks=true, pdfstartview=FitV, linkcolor=blue, citecolor=blue, urlcolor=blue,bookmarks,bookmarksnumbered]{hyperref}
\hypersetup{
pdftitle={Laplace},
pdfauthor={S. Gill Williamson}%{\textcopyright\ S. Gill Williamson}
}
\begin{document}
\thispagestyle{empty}
\begin{center}
\vspace*{1in}
\textsf{ \Large The common-submatrix Laplace expansion}\\[.2in]
%\textsf{\small basic skills and proof techniques}\\[.7in]

\textsf{S. Gill Williamson}
%\vfill
%\textcopyright  \textsf{S. Gill Williamson 2012. All rights reserved.}
\end{center}

%\newpageI
\thispagestyle{empty}
\hspace{1 pt}
%\newpage
%ABSTRACT
\begin{center}
{\Large Abstract}\\[.2in]
\end{center}
\pagestyle{plain}
We state and prove a classical version of  the Laplace expansion theorem where  all submatrices in the expansion are restricted to contain a specified {\em common submatrix (CSM)}. 
The result states that (accounting  for signs) this restricted expansion equals the determinant of the original matrix times the determinant of the CSM.  
This result (Muir~\cite{tm:ttd}, p.132) is one of many such results contained in {\em A Treatise on the Theory of Determinants} by Muir (revised and enlarged by Metzler).  
Our proof,  based on the same general idea used by Muir and valid for any commutative ring, is a modification of  a proof used in Williamson~\cite{gw:gcm} to study the non-commuting case.   
 
%(\cite{gw:gcm}, p. 254, 
%INTRODUCTION
\section{\bfseries Introduction}
The book~\cite{tm:ttd} by Muir and Metzler contains an amazing collection of combinatorial methods applied to determinantal identities.  All of Muir's results and proofs are stated in standard English and proved by example, making it difficult to decipher statements and verify proofs. 
In this paper, we state and prove the CSM Laplace expansion, 
Muir~\cite{tm:ttd} (p. 132) and theorem~\ref{thm:csmlap} below, using the standard Laplace expansion (theorem~\ref{thm:rnkvrsstn}). 
Our definitions, statements and proofs use elementary set theoretic and combinatorial methods and, thus, give a pattern for the translation of other of  Muir's results~\cite{tm:ttd} into more familiar mathematics.
We state our results for matrices where the rows and columns are indexed by arbitrary linearly ordered sets.  This level of generalization is useful for certain applications and in our proofs.

We need some notational conventions.
%REM
\begin{remark}[\bfseries Notation]
$\bP_k(n)$ is the set of all subsets $X$  of $\underline{n}\coloneq\{1, \ldots ,n\}$ of size $k$ (i.e., $|X|=k$).
We write $A\in\mathbf{M}_{\underline{n}}$ or  $A\in\mathbf{M}_{n}$
to designate an $n\times n$ matrix with entries in $\bK$, a commutative ring.
 We use $A[X\vert Y]$ to denote the submatrix of $A$ gotten by retaining rows indexed
 by the set $X$ and columns indexed by the set $Y$.
 We use $A(X\vert Y)$ to denote the submatrix of $A$ gotten by retaining rows indexed by the set $\underline{n}\setminus X$ (the complement of $X$ in $\underline{n}$) and columns indexed by the set 
 $\underline{n}\setminus Y$.  
 We also use the mixed notation  $A[X\vert Y)$ and  $A(X\vert Y]$
 with obvious meaning.
 We use $\Theta$ to denote the zero matrix of the appropriate size and $I$ to denote the identity matrix.
 We abbreviate $\sum_{x\in S} f(x)$ to $\sum_{S} f(x)$.
\end{remark}
%END REM
\begin{definition}[\bfseries Position and rank functions for linear orders]
\label{rem:pstrnk}
Let $(\Lambda, \leq)$ be a finite linearly ordered set.
For $x\in \Lambda$ and $S\subseteq \Lambda$ let
\begin{equation}
\pi^{\Lambda}_{S}(x)\coloneq\left|\{t\Mid t\in S, t\leq x\}\right|\;\;\mathrm{and}\;\;
\rho^{\Lambda}_{S}(x)\coloneq\left|\{t\Mid t\in S, t < x\}\right|.
\end{equation}
$\pi^{\Lambda}_{S}$ is called the {\em position function} for $\Lambda$ relative  to $S$ and 
$\rho^{\Lambda}_S$ the {\em rank function} for $\Lambda$ relative to $S$.
If $S=\Lambda$, we use $\pi^{\Lambda}$ instead of $\pi^{\Lambda}_{\Lambda}$, and,
similarly, we use $\rho^{\Lambda}$ instead of $\rho^{\Lambda}_{\Lambda}$.
\end{definition}
The standard Laplace expansion, general indices, is as follows:
%THM
\begin{theorem}[\bfseries Standard Laplace expansion, general indices case]
\label{thm:rnkvrsstn}
Let $A$ be a matrix with entries in a commutative ring, rows and columns 
labeled by linearly ordered sets  $\Phi$ and $\Lambda$ respectively, 
$|\Phi|=|\Lambda|=n$.
Let $(K, K')$ and $(L, L')$ be  ordered set partitions of $\Phi$ and $\Lambda$ where $|K|=|L|$, $|K'|=|L'|$.  
We assume that $K$ (hence $K'$) is fixed and $L$ (hence $L'$) is variable. 
We have $\;\det(A) = $
\begin{equation}
\label{eq:anc1}
(-1)^{\sum_{K}\pi^{\Phi}(x)} \sum_{L} (-1)^{\sum_L\pi^{\Lambda}(x)}\det(A[K\Mid L])\det(A[K'\vert L'])
\end{equation}
or, alternatively,
\begin{equation}
\label{eq:anc2}
(-1)^{\sum_{K}\rho^{\Phi}(x)} \sum_{L} (-1)^{\sum_L\rho^{\Lambda}(x)}
\det(A[K\Mid L])\det(A[K'\vert L']).
\end{equation}
One can replace $\det(A[K'\vert L'])$ with $\det(A(K\Mid L))$ in these identities.
\end{theorem}
Theorem~\ref{thm:rnkvrsstn}, is a classical result in the case $\Phi=\Lambda=\underline{n}$ which
is proved in many texts.   For an online proof source, see~\cite{gw:mcf}, p $49$.
The general indices case follows from the natural order isomorphisms of ${\Phi}$ with 
$\{\pi^{\Phi}(x)\Mid x\in \Phi\}$ and  ${\Lambda}$ with $\{\pi^{\Lambda}(x)\Mid x\in \Lambda\}.$ 
Equation~\ref{eq:anc2} follows from 
\[
\sum_K\pi^{\Phi}(x) = \sum_K\rho^{\Phi}(x) + |K|,\;
\sum_L\pi^{\Lambda}(x) = \sum_L\rho^{\Lambda}(x) + |L|,\; \mathrm{and}\; |K|=|L|.
\]
 The following generalization of theorem~\ref{thm:rnkvrsstn} will be the focus of this paper.
%anchored laplace    
\begin{theorem}[\bfseries CSM Laplace expansion, general indices case]
\label{thm:csmlap}
Let $A$ be a matrix with entries in a commutative ring, rows and columns 
labeled by linearly ordered sets  $\Phi$ and $\Lambda$ respectively, $|\Phi|=|\Lambda|=n$. 
Let $(F, I, I')$ and $(G, J, J')$ be  ordered set partitions of  $\Phi$ and $\Lambda$ respectively, where
$|F|=|G|$,  $|I|=|J|$ (hence $|I'|=|J'|$). We assume $F$, $G$, $I$ (hence $I'$) are fixed and
$J$ (hence $J'$) is variable. 
We have 
$\;\det(A[F\mid G])\det(A) =$
\begin{equation}
\label{eq:csmlap1}
(-1)^{\sum_{I}\pi^{\Phi}_{I\cup I'}(x)}\sum_{J} (-1)^{\sum_{J}\pi^{\Lambda}_{J\cup J'}(x)}
 \det(A[F\cup I\Mid G\cup J])\det A[F\cup I'\Mid G\cup J']
\end{equation}
or, alternatively, 
\begin{equation}
\label{eq:csmlap2}
(-1)^{\sum_{I}\rho^{\Phi}_{I\cup I'}(x)}\sum_{J} (-1)^{\sum_{J}\rho^{\Lambda}_{J\cup J'}(x)}
 \det(A[F\cup I\mid G\cup J])\det A[F\cup I'\Mid G\cup J']. 
\end{equation}
We refer to the matrix $A[F\,|\,G]$ as the {\em common submatrix} of the expansion, noting
that it is a submatrix of both $A[F\cup I\mid G\cup J]$ and  
$A[F\cup I'\,|\, G\cup J']=A(I\,|\, J)$.
\begin{proof}
The case $\Phi = \Lambda = \underline{n}$ is  lemma~\ref{lem:prfcsmlpl}.  
The case for general $\Phi$ and $\Lambda$ is routine and follows from the natural order isomorphisms of ${\Phi}$ with  $\{\pi^{\Phi}(x)\Mid x\in \Phi\}$ and  ${\Lambda}$ with $\{\pi^{\Lambda}(x)\Mid x\in \Lambda\}.$  
\end{proof}  
\end{theorem}
%SECTION
\section{\bfseries Statements and proofs}
As in the hypothesis of theorem~\ref{thm:csmlap}, let $A$ be an $n\times n$ matrix with entries in a commutative ring.
Let $(F, I, I')$ and $(G, J, J')$ be  ordered set partitions of  $\underline{n}$ where
$|F|=|G|$,  $|I|=|J|$ (hence $|I'|=|J'|$). We assume $F$, $G$, $I$ (hence $I'$) are fixed and
$J$ (hence $J'$) is variable.  
%DEF
\begin{definition}[\bfseries Extending the natural order on $\underline{n}$]
\label{def:sbsxtnlnr}
Let $F\subseteq \underline{n}$.  Define $F^{+}=\{f^{+}\mid f\in F\}$ where the notation $f^{+}$ represents a new symbol that is the successor to $f$ in a linear  order on $\underline{n}^{+}_F \coloneq \underline{n}\cup F^{+}$ which restricts to the natural order on $\underline{n}$.
Thus, for $n=6$ and $F=\{2, 4\}$, 
$\underline{n}^{+}_F = (1,2, 2^{+}, 3, 4, 4^{+}, 5, 6)$.
We use $\mathbf{M}_{\underline{n}^{+}_F,\underline{n}^{+}_G}$ to denote the matrices with rows indexed by 
${\underline{n}^{+}_F}$ and columns indexed 
by ${\underline{n}^{+}_G}$, $G\subseteq \underline{n}$, $|G|=|F|$  (see~\ref{eq:mtrx02}).
\end{definition}

\begin{remark}[\bfseries Basic facts about rank functions and $\underline{n}^{+}_F$]
\label{rem:bsfcts}
%EQ
Equation~\ref{eq:bsfcts1} holds for any disjoint union $S_1\dot\cup S_2={\underline{n}^{+}_F}$, in particular for $S_1=F\cup F^{+}$,
$S_2=I\cup I'$.
\begin{equation}
\label{eq:bsfcts1}
(F\cup F^{+}) \dot\cup (I\cup I') = {\underline{n}^{+}_F}\;\mathrm{and}\; x\in\underline{n}^{+}_F\implies 
\rho^{\underline{n}^{+}_F}(x)= 
 \rho^{\underline{n}^{+}_F}_{F\cup F^{+}}(x) + \rho^{\underline{n}^{+}_F}_{I\cup I'}(x)
 \end{equation}
  %EQ
 Equation~\ref{eq:bsfcts2} follows directly from definition~\ref{rem:pstrnk} and~\ref{def:sbsxtnlnr}. 
\begin{equation}
\label{eq:bsfcts2}
S\subseteq \underline{n},\; x\in\underline{n} \implies \rho_S^{\underline{n}^{+}_F}(x) =\rho_S^{\underline{n}}(x)
\end{equation}
  %EQ
 Equation~\ref{eq:bsfcts3} follows from the adjacent pairing of the $f_if_i^{+}$. 
\begin{equation}
\label{eq:bsfcts3}
x\in\underline{n} \implies \rho_{F\cup F^{+}}^{\underline{n}^{+}_F}(x) = 0\, (\mathrm{mod}\; 2)
\end{equation}
 %EQ
 Equation~\ref{eq:bsfcts4} follows from equation~\ref{eq:bsfcts1}, \ref{eq:bsfcts2} and \ref{eq:bsfcts3}
\begin{equation}
\label{eq:bsfcts4}
x\in\underline{n} \implies
\rho^{\underline{n}^{+}_F}(x)=\rho_{I\cup I'}^{\underline{n}^{+}_F}(x)= \, 
\rho_{I\cup I'}^{\underline{n}}(x)\,(\mathrm{mod}\; 2)
\end{equation}
\end{remark}

\begin{sremark}[\bfseries Example of {``initializing''} matrices]
\label{rem:trnsfmts}
We use the notation of theorem~\ref{thm:csmlap}. We start with $A\in \mathbf{M_6}$ and transform it to a larger {\em initialized matrix} $\vec{A}$ (definition~\ref{def:intlmtr} specifies $\vec{A}$ in general).  
We first transform $A$ to $\widehat{A}$ by symbolically ``doubling'' the rows indexed by $F$ and the columns indexed by $G$.
Here $F=\{f_1,f_2\}=\{2,4\}$, $F^{+}\coloneq\{f_1^{+}, f_2^{+}\}=\{2^{+},4^{+}\}$, $G=\{g_1,g_2\}=\{3,5\}$, 
$G^{+}\coloneq \{g_1^{+}, g_2^{+}\}=\{3^{+},5^{+}\}$, $I=\{i_1,i_2\}=\{1,6\}$, and $I'=\{i'_1,i'_2\}=\{3,5\}$.
The extended linear order for rows is $\underline{n}_F^{+}=(1,2,2^{+},3,4,4^{+},5,6)$ (see~\ref{def:sbsxtnlnr}).
We replace $A(i,j)$ by $\mathrm{ij}$ for convenience of notation.
\begin{equation}
\label{eq:mtrx01}
A =\, 
\bordermatrix{~ &~ & ~ & \overset{g_1}{3}& ~&\overset{g_2}{5}&~ \cr
                  i_1=1 & 11 & 12 & 13 & 14 & 15 & 16\cr
                  f_1=2& 21 & 22 & 23 &24 & 25 & 26\cr
                  i'_1=3 & 31 & 32 & 33 & 34 & 35 & 36\cr
                  f_2=4& 41 & 42 & 43 & 44 & 45 & 46\cr
                  i'_2=5& 51 & 52 & 53 & 54 & 55 & 56\cr
                  i_2=6 & 61 & 62 & 63 & 64 & 65 & 66\cr}
\end{equation}
\begin{equation}
\label{eq:mtrx02}
\widehat{A} =\, 
\bordermatrix{~ &~ & ~ & 3 & 3^{+}&~ & 5 & 5^{+} &~ \cr
                i_1= 1   & 11 & 12 & 13 & 13 & 14 & 15 & 15 & 16\cr
                 f_1= 2   & 21 & 22 & 23 & 23 & 24 & 25 & 25 & 26\cr
                  f_1^{+}=2^{+} & 21 & 22 & 23 & 23 & 24 & 25 & 25 & 26\cr
                  i'_1= 3   & 31 & 32 & 33 & 33 & 34 & 35 & 35 & 36\cr
                  f_2=4   & 41 & 42 & 43 & 43 & 44 & 45 & 45 & 46\cr
                  f_2^{+}=4^{+} & 41 & 42 & 43 & 43 & 44 & 45 & 45 & 46\cr
                 i'_2= 5   & 51 & 52 & 53 & 53 & 54 & 55 & 55 & 56\cr
                  i_2=6   & 61 & 62 & 63 & 63 & 64 & 65 & 65 & 66\cr}
\end{equation}
\begin{equation}
\label{eq:mtrx03}
 \vec{A} =\, 
\bordermatrix{~ &~ & ~ & 3 & 3^{+}&~ & 5 & 5^{+} &~ \cr
                  i_1= 1   & 11 & 12 & 13 & 0 & 14 & 15 & 0 & 16\cr
                  f_1= 2   & 21 & 22 & 23 & 0 & 24 & 25 & 0 & 26\cr
                  f_1^{+}=2^{+}  & 21 & 22 & 0 & 23 & 24 & 0 & 25 & 26\cr
                  i'_1= 3   & 31 & 32 & 0 & 33 & 34 & 0 & 35 & 36\cr
                  f_2=4   & 41 & 42 & 43 & 0 & 44 & 45 & 0 & 46\cr
                  f_2^{+}=4^{+} & 41 & 42 & 0 & 43 & 44 & 0 & 45 & 46\cr
                  i'_2= 5  & 51 & 52 & 0 & 53 & 54 & 0 & 55 & 56\cr
                   i_2=6  & 61 & 62 & 63 & 0 & 64 & 65 & 0 & 66\cr}
\end{equation}
\end{sremark}
%END SREM
%SREM
\begin{srem}[\bfseries Example of equivalent transformations]
\label{srem:dtmeqvtrn}
From (~\ref{eq:mtrx03}), let  
\begin{equation}
\label{eq:mtrx04}
\vec{A} =\, 
\bordermatrix{~ &1 &   2  & 3 &  3^{+}&  4&    5 & 5^{+} & 6 \cr
                  1   & 11 & 12 & 13 & 0 & 14 & 15 & 0 & 16\cr
                  2   & 21 & 22 & 23 & 0 & 24 & 25 & 0 & 26\cr
                  2^{+} & 21 & 22 & 0 & 23 & 24 & 0 & 25 & 26\cr
                  3   & 31 & 32 & 0 & 33 & 34 & 0 & 35 & 36\cr
                  4   & 41 & 42 & 43 & 0 & 44 & 45 & 0 & 46\cr
                  4^{+} & 41 & 42 & 0 & 43 & 44 & 0 & 45 & 46\cr
                  5   & 51 & 52 & 0 & 53 & 54 & 0 & 55 & 56\cr
                  6   & 61 & 62 & 63 & 0 & 64 & 65 & 0 & 66\cr}
\end{equation}
We apply elementary column and row operations to $\vec{A}$. 
%(\cite{gw:mcf}, p $63$, definition~4.7). 
Let $\hat{R}_{[i]-[j]}X$ denote replacing row $i$ of $X$ by the difference of row $i$ and $j$.
Let $\hat{C}_{[i]+[j]}X$ denote replacing column $i$ of $X$ by the sum of column 
$i$ and $j$.
Let $F= \{2, 4\}$, $F^{+}=\{2^{+}, 4^{+}\}$, $I=\{1, 6\}$ and $I'=\{3,5\}$.
Let $G=\{3, 5\}$ and $G^{+} = \{3^{+}, 5^{+}\}$. 
Let $\hat{C}_{[G]+[G^{+}]} =\{\hat{C}_{[i]+[i^{+}]}\Mid i\in G\}$ (a set of commuting operators).
By $\hat{C}_{[G]+[G^{+}]}\vec{A}$ we mean applying these operators to $\vec{A}$ in any order.
\begin{equation}
\label{eq:mtrx05}
\hat{C}_{[G]+[G^{+}]}\vec{A} =\, 
\bordermatrix{~ &1 &   2  & 3 &  3^{+}&  4&    5 & 5^{+} & 6 \cr
                  1   & 11 & 12 & 13 & 0 & 14 & 15 & 0 & 16\cr
                  2   & 21 & 22 & 23 & 0 & 24 & 25 & 0 & 26\cr
                  2^{+} & 21 & 22 & 23 & 23 & 24 & 25 & 25 & 26\cr
                  3   & 31 & 32 & 33 & 33 & 34 & 35 & 35 & 36\cr
                  4   & 41 & 42 & 43 & 0 & 44 & 45 & 0 & 46\cr
                  4^{+} & 41 & 42 & 43 & 43 & 44 & 45 & 45 & 46\cr
                  5   & 51 & 52 & 53 & 53 & 54 & 55 & 55 & 56\cr
                  6   & 61 & 62 & 63 & 0 & 64 & 65 & 0 & 66\cr}
\end{equation}
and, setting $\widetilde{A} = \hat{R}_{[F^{+}]-[F]}\hat{C}_{[G]+[G^{+}]}\vec{A}$, we have
\begin{equation}
\label{eq:mtrx06}
\widetilde{A}  =\, 
\bordermatrix{~ &1 &   2  & 3 &  3^{+}&  4&    5 & 5^{+} & 6 \cr
                  1   & 11 & 12 & 13 & 0 & 14 & 15 & 0 & 16\cr
                  2   & 21 & 22 & 23 & 0 & 24 & 25 & 0 & 26\cr
                  2^{+} & 0   & 0   & 0   & 23 & 0 & 0 & 25 & 0\cr
                  3   & 31 & 32 & 33 & 33 & 34 & 35 & 35 & 36\cr
                  4   & 41 & 42 & 43 & 0 & 44 & 45 & 0 & 46\cr
                  4^{+} & 0   & 0   & 0   &43& 0   & 0 & 45 & 0\cr
                  5   & 51 & 52 & 53 & 53 & 54 & 55 & 55 & 56\cr
                  6   & 61 & 62 & 63 & 0 & 64 & 65 & 0 & 66\cr}.
\end{equation}
Note that 
\[ 
|\det(\widetilde{A})|=|\det (\widetilde{A}[F^{+}\mid G^{+}])\det(\widetilde{A}(F^{+}\mid G^{+})| =
|\det(A[F\mid G]) \det(A)|
\]
which is the left hand side of~\ref{eq:csmlap1} up to sign. 
Note that $\det(\widetilde{A}) = \det(\vec{A})$ since one can be transformed into the other by ``type II'' elementary row and column operations which preserve determinants.
We shall show that $\det(\vec{A})$ yields the right hand side  of~\ref{eq:csmlap1}. 
\end{srem}
%END SREM
%DEF
\begin{definition}[\bfseries Initialization matrix $\vec{A}$ of $A$]
\label{def:intlmtr}
Given $A\in \mathbf{M}_{\underline{n}}$, we define 
$\vec{A}\in\mathbf{M}_{\underline{n}^{+}_F,\underline{n}^{+}_G}$ by
 the following conditions:
\begin{equation}
\label{eq:intlmtr1}
%\vec{A}[F\cup I\Mid G] = {A}[F\cup I\Mid G] \;\;\mathrm{and}\;\; 
\vec{A}[F\cup I\Mid G^{+}) = {A}[F\cup I\Mid \underline{n}] \;\mathrm{and}\;
\vec{A}[F\cup I\Mid G^{+}] =\Theta,\;
\end{equation}
\begin{equation}
\label{eq:intlmtr2}
\vec{A}[F^{+}\cup I'\Mid G) = {A}[F\cup I'\Mid \underline{n}] \;\mathrm{and}\;
\vec{A}[F^{+}\cup I'\Mid G] =\Theta,\;
\end{equation}
These conditions are easily checked in the example given in~\ref{eq:mtrx03}. 
\end{definition}
%REM
\begin{remark}[\bfseries Observations about $\vec{A}$]
\label{rem:obsrvtns}
With the notation of~\ref{thm:csmlap}, 
$(F, I, I')$ and $(G, J, J')$, ordered partitions of  $\underline{n}$, become,
using definition~\ref{def:sbsxtnlnr}, 
$(F, F^{+},I, I')$ and $(G,G^{+}, J, J')$,
ordered partitions of $\Phi=\underline{n}^{+}_F$ and 
$\Lambda=\underline{n}^{+}_G$ respectively where
$|F|=|F^{+}|=|G|=|G^{+}|$, $|I|=|J|$, $|I'|=|J'|$.
We refer to definition~\ref{def:intlmtr} and the example shown in equation~\ref{eq:mtrx03}.
It is easily seen, that   
the conditions of~\ref{eq:intlmtr1} and~\ref{eq:intlmtr2}  imply that
\begin{equation}
\label{eq:obsrvtns1}
\vec{A}[F\cup I \vert G\cup J] = A[F\cup I \vert G\cup J]
\end{equation}
%EQ
\begin{equation}
\label{eq:obsrvtns2}
\vec{A}[F^{+}\cup I' \vert G^{+}\cup J'] = A[F\cup I' \vert  G\cup J'].
\end{equation}
These relations are evident for $\vec{A}$ of~\ref{eq:mtrx03}. 
Let $(L,L')$ be a partition of  $\Lambda$ where $|L|=|F\cup I|$.
From the definition of $\vec{A}$ (definition~\ref{def:intlmtr}), note that
 if $L\cap G^{+} \neq \emptyset$ then $\vec{A}[F\cup I \mid L]$ contains a zero column.  
Similarly, $L'\cap G \neq \emptyset$ implies $\vec{A}[F^{+}\cup I' \mid L'])$
contains a zero column.
Thus, 
\begin{equation}
\label{eq:obsrvtns3}
\det(\vec{A}[F\cup I \Mid L])\det(\vec{A}[F^{+}\cup I' \Mid L'])\neq 0 \implies
(L,L')=(G\cup J, G^{+}\cup J')
\end{equation}
where $J=L\setminus G$ and $J'=L'\setminus G^{+}$.
This $(G,G^{+}, J, J')$
defines a unique ordered partition of 
$\Lambda=\underline{n}^{+}_G$ where
$|F|=|F^{+}|=|G|=|G^{+}|$, $|I|=|J|$, $|I'|=|J'|$.
\end{remark}
%REM

\begin{lemma}[\bfseries Right hand side of CSM Laplace]
\label{lem:appstnlpl}
We use the notation of remark~\ref{rem:obsrvtns}.
Let $\vec{A}$ be as in definition~\ref{def:intlmtr}.
Then $\det(\vec{A}) =$
\begin{equation*}
(-1)^{\sum_{F\cup I} \rho^{\underline{n}}_{I\cup I'}(x)} \sum_{J} (-1)^{\sum_{G\cup J}\rho^{\underline{n}}_{J\cup J'}(x)}\det(A[F\cup I\mid G\cup J])\det(A[F\cup I'\Mid G\cup J']).
\end{equation*}
\begin{proof}
Recall that $(F, F^{+}, I, I')$, $G$ and $G^{+}$ are fixed and $J$ and $J'$ range over all ordered 
partitions of the form $(G, G^{+}, J, J')$ where $|J|=|I|$ and thus  $|J'|=|I'|$.
Let  $\Phi=\underline{n}^{+}_F$, 
$\Lambda=\underline{n}^{+}_G$ be as in definition~\ref{def:sbsxtnlnr}.
Let $(K,K')$ and $(L,L')$ be as in theorem~\ref{thm:rnkvrsstn}. Take $K=F\cup I$, $K'=F^{+}\cup I'$. 
From~\ref{thm:rnkvrsstn} we have
\begin{equation}
\label{eq:appstnlpl}
\det(\vec{A}) = (-1)^{\sum_{K}\rho^{\Phi}(x)} \sum_{L\subseteq \Lambda} (-1)^{\sum_L\rho^{\Lambda}(x)}
\det(\vec{A}[K\Mid L])\det(\vec{A}[K'\Mid L']).
\end{equation}
From equation~\ref{eq:obsrvtns3} of remark~\ref{rem:obsrvtns}, equation~\ref{eq:appstnlpl} becomes
\begin{equation}
\label{eq:appstnlpl0} 
\det(\vec{A}) =
f(x)\sum_{J} g_J(x) \det(\vec{A}[F\cup I\mid G\cup J])
\det(\vec{A}[F^{+}\cup I'\Mid G^{+}\cup J'])
\end{equation}
where $f(x)=(-1)^{\sum_{F\cup I}\rho^{\Phi}(x)}$ and $g_J(x)=(-1)^{\sum_{G\cup J} \rho^{\Lambda}(x)}.$
From remark~\ref{rem:obsrvtns}, we now have
\begin{equation}
\label{eq:appstnlpl1} 
\det(\vec{A}) =
f(x)\sum_{J} g_J(x) \det(A[F\cup I\mid G\cup J])\det(A[F\cup I'\Mid G\cup J']).
\end{equation}
Noting that $F\cup I\subseteq \underline{n}$ and using~\ref{rem:bsfcts} (equation~\ref{eq:bsfcts4}),
$f(x)=(-1)^{\sum_{F\cup I}\rho_{I\cup I'}^{\underline{n}}}(x)$ and, similarly,  
$g_J(x)=(-1)^{\sum_{G\cup J}\rho_{J\cup J'}^{\underline{n}}}(x).$                 
\end{proof}
\end{lemma}

%LEM
\begin{lemma}[\bfseries Left hand side of CSM Laplace]
\label{lem:lhscsm}
Referring to the notation of lemma~\ref{lem:appstnlpl}, the following identity holds.
\begin{equation}
\label{eq:lhscsm1}
\det(\vec{A})= (-1)^{\sum_F\rho^{\underline{n}}_{I\cup I'}(x)} (-1)^{\sum_G\rho^{\underline{n}}_{J\cup J'}(x)}
\det(A[F\Mid G]\det(A).
\end{equation}
\begin{proof}
Let  $\Phi=\underline{n}^{+}_F$, 
$\Lambda=\underline{n}^{+}_G$ be as in definition~\ref{def:sbsxtnlnr}.
The proof follows example~\ref{srem:dtmeqvtrn}.
\begin{equation}
\label{eq:lhscsm1}
\widetilde{A} = \hat{R}_{[F^{+}]-[F]}\hat{C}_{[G]+[G^{+}]}\vec{A}.
\end{equation}
Taking $K=F^{+}$ in theorem~\ref{thm:rnkvrsstn}, we find only one term in the expansion:
\begin{equation}
\label{eq:lhscsm2}
\det(\widetilde{A})=(-1)^{\sum_{F^{+}}\rho^{\Phi}(x)} (-1)^{\sum_{G^{+}}\rho^{\Lambda}(x)}
\det(\widetilde{A}[F^{+}\Mid G^{+}])\det(\widetilde{A}(F^{+}\Mid G^{+})).
\end{equation}
By example~\ref{srem:dtmeqvtrn}, the general case follows the same pattern: 
$\widetilde{A}[F^{+}\Mid G^{+}]=A[F\Mid G]$ and $\widetilde{A}(F^{+}\Mid G^{+})=A.$
Note that $\sum_{F^{+}}\rho^{\Phi}(x)=|F|+ \sum_{F}\rho^{\Phi}(x)$ and
$\sum_{G^{+}}\rho^{\Phi}(x)=|G|+ \sum_{G}\rho^{\Phi}(x).$ 
Since $|F|=|G|$ we have 
\begin{equation}
\label{eq:lhscsm3}
\det(\widetilde{A})=(-1)^{\sum_{F}\rho^{\Phi}(x)} (-1)^{\sum_{G}\rho^{\Lambda}(x)}\det(A[F\Mid G])\det(A).
\end{equation} 
Note that $F\subseteq \underline{n}.$  Using~\ref{rem:bsfcts} (equation~\ref{eq:bsfcts4}),
$(-1)^{\sum_{F}\rho^{\Phi}(x)} =(-1)^{\sum_{F}\rho_{I\cup I'}^{\underline{n}}}(x)$ and, similarly,  
$(-1)^{\sum_{G}\rho^{\Lambda}(x)}=(-1)^{\sum_{G}\rho_{J\cup J'}^{\underline{n}}}(x).$ 
\end{proof}
\end{lemma}
%LEM
\begin{lemma}[\bfseries CSM Laplace for $\underline{n}$]
\label{lem:prfcsmlpl}
We have 
$\;\det(A[F\mid G])\det(A) =$
\begin{equation*}
(-1)^{\sum_{I}\rho^{\underline{n}}_{I\cup I'}(x)}\sum_{J} (-1)^{\sum_{J}\rho^{\underline{n}}_{J\cup J'}(x)}
 \det(A[F\cup I\Mid G\cup J])\det (A[F\cup I' \Mid G\cup J'])
\end{equation*}
or, alternatively,
\begin{equation*}
(-1)^{\sum_{I}\rho^{\underline{n}}_{I\cup I'}(x)}\sum_{J} (-1)^{\sum_{J}\rho^{\underline{n}}_{J\cup J'}(x)}
 \det(A(I'| J'))\det (A(I\,|\, J)). 
\end{equation*}
\begin{proof}
Consider the  two expressions for $\det(\vec{A})$ from lemma~\ref{lem:appstnlpl} and lemma~\ref{lem:lhscsm}
which are, respectively,
\begin{equation*}
(-1)^{\sum_{F\cup I} \rho^{\underline{n}}_{I\cup I'}(x)} \sum_{J} (-1)^{\sum_{G\cup J}\rho^{\underline{n}}_{J\cup J'}(x)}\det(A[F\cup I\mid G\cup J])\det(A[F\cup I'\Mid G\cup J']
\end{equation*}
and
\begin{equation*}
(-1)^{\sum_F\rho^{\underline{n}}_{I\cup I'}(x)} (-1)^{\sum_G\rho^{\underline{n}}_{J\cup J'}(x)}
\det(A[F\Mid G]\det(A).
\end{equation*}
The sums from lemma~\ref{lem:appstnlpl}  satisfy
\[
{\sum_{F\cup I} \rho^{\underline{n}}_{I\cup I'}(x)}=\sum_{F} \rho^{\underline{n}}_{I\cup I'}(x) +
\sum_{I} \rho^{\underline{n}}_{I\cup I'}(x)
\]
\[
{\sum_{G\cup J} \rho^{\underline{n}}_{J\cup J'}(x)}=\sum_{G} \rho^{\underline{n}}_{J\cup J'}(x) +
\sum_{J} \rho^{\underline{n}}_{J\cup J'}(x).
\]
The sums $\sum_{F}\rho^{\underline{n}}_{I\cup I'}(x)$ and 
$\sum_{G} \rho^{\underline{n}}_{J\cup J'}(x)$ occur in the expression in lemma~\ref{lem:lhscsm}.
Canceling appropriate terms from the expressions for $\det(\vec{A})$ from lemma~\ref{lem:appstnlpl} and lemma~\ref{lem:lhscsm} completes the proof.
\end{proof}
\end{lemma}
%SECTION
\section{Examples}
\begin{sremark}[\bfseries Example]
\label{srem:usgprf}
Consider the  transformations: $A \rightarrow \widehat{A} \rightarrow \vec{A}$ together with 
$\widetilde{A} =  \hat{R}_{[F^{+}]-[F]}\hat{C}_{[G]+[G^{+}]}\vec{A}.$ 
 (\ref{eq:mtrx01}, \ref{eq:mtrx02}, \ref{eq:mtrx03} and \ref{eq:mtrx06} ).

For $F=\{1\},\,G=\{1\}\;$:
\[
A=(a_{11}) 
\rightarrow 
\widehat{A}=
\bordermatrix
{
~      &   1        &1^{+}\cr
1      &a_{11}   &a_{11}\cr
1^{+}&a_{11}   &a_{11}\cr
}
\;\;
\vec{A}=
\bordermatrix
{
~      &   1        &1^{+}\cr
1      &a_{11}   &0\cr
1^{+}&0           &a_{11}\cr
}
\]
\[
\widetilde{A} = \hat{R}_{[F^{+}]-[F]}\hat{C}_{[G]+[G^{+}]}\vec{A}=
\bordermatrix
{
~      &   1        &1^{+}\cr
1      &a_{11}   &0\cr
1^{+}&0           &a_{11}\cr
}.
\]

\end{sremark}
%REM
\begin{sremark}[\bfseries Example]
\label{srem:usgprf}
Consider CSM Laplace for $\underline{n}$ with $n=3$ (lemma~\ref{lem:prfcsmlpl}).
Let $I=\{1\},\;F=\{2\},\;I'=\{3\}\;$ and $G=\{2\}.$
\[
A=
\bordermatrix
{
~      &   1        &2  &                3\cr
1      &a_{11}   &a_{12}  &a_{13}\cr
2      &a_{21}   &a_{22}  &a_{23}\cr
3      &a_{31}   &a_{32}  &a_{33}\cr
}
\]
The RHS of lemma~\ref{lem:prfcsmlpl} has two terms, for $J=\{1\}$ and $\{J\}=\{3\}$:
\begin{equation*}
(-1)^{\sum_{I}\rho^{\underline{n}}_{I\cup I'}(x)}\sum_{J} (-1)^{\sum_{J}\rho^{\underline{n}}_{J\cup J'}(x)}
 \det(A[F\cup I\mid G\cup J])\det (A[F\cup I'\,|\, G\cup J'])=
 \end{equation*}
 \[
(-1)^{\sum_{\{1\}}\rho^{\underline{3}}_{\{1,3\}}(1)}
\left((-1)^{\sum_{\{1\}}\rho^{\underline{3}}_{\{1,3\}}(1)}\det(A[\{1,2\}\mid \{1,2\}])\det (A[\{2,3\}\,|\, \{2,3\}])\right)+
\]
\[
(-1)^{\sum_{\{1\}}\rho^{\underline{3}}_{\{1,3\}}(1)}
\left((-1)^{\sum_{\{3\}}\rho^{\underline{3}}_{\{1,3\}}(3)}\det(A[\{1,2\}\mid \{2,3\}])\det (A[\{2,3\}\,|\, \{1,2\}])\right)=
\]
\[
(-1)^{0}
\left((-1)^{0}\det(A[\{1,2\}\mid \{1,2\}])\det (A[\{2,3\}\,|\, \{2,3\}])\right)+
\]
\[
(-1)^{0}
\left((-1)^{1}\det(A[\{1,2\}\mid \{2,3\}])\det (A[\{2,3\}\,|\, \{1,2\}])\right)=
\]
\[
(a_{11} a_{22} - a_{21} a_{12})(a_{22} a_{33} - a_{32} a_{23}) - 
(a_{12} a_{23} - a_{22} a_{13})(a_{21} a_{32} - a_{31} a_{22}).
\]
Note that $\det (A[F\Mid G]) = a_{22}$, thus, from  lemma~\ref{lem:prfcsmlpl},  we have
\[
a_{22}\det(A) = (a_{11} a_{22} - a_{21} a_{12})(a_{22} a_{33} - a_{32} a_{23}) - 
(a_{12} a_{23} - a_{22} a_{13})(a_{21} a_{32} - a_{31} a_{22}).
\]
Expanding the RHS of the above, we get eight terms:
\[
a_{11}\underline{a_{22}}a_{22} a_{33} - a_{11}\underline{a_{22}}a_{32} a_{23} - a_{21} a_{12}\underline{a_{22}} a_{33} + 
\overbrace{a_{21} a_{12}a_{32} a_{23}}
\]
\[
-\overbrace{a_{12} a_{23}a_{21} a_{32}} + a_{12} a_{23}a_{31}\underline{a_{22}} + \underline{a_{22}} a_{13}a_{21} a_{32} - 
\underline{a_{22}} a_{13}a_{31} a_{22}.
\]
Canceling the terms with over-braces and factoring out the $\underline{a_{22}}$ terms gives
$a_{22}\det(A)$ as required.
\end{sremark}
%SREM
\begin{sremark}[\bfseries Example]
\label{srem:fnlxmp}
Example~\ref{srem:usgprf} generalizes in an interesting way.
Consider CSM Laplace for $\underline{n}$ (lemma~\ref{lem:prfcsmlpl}).
Let $I=\{1\}$,$\;F=\{2,\ldots, n-1\}$, $\;I'=\{n\}\;$ and $\;G=\{2,\ldots, n-1\}.$
The pairs $(J, J') \in \{(\{1\}, \{n\}),\; (\{n\}, \{1\})\}$ have just two possible choices.
We set $A[F\Mid G]$ as $A(\{1,n\}\Mid \{1,n\}) \coloneq A(1,n\Mid 1,n)$, etc.
Lemma~\ref{lem:prfcsmlpl} (second identity) becomes
\begin{equation*}
\label{eq:fnlxmp1}
\det(A(1,n\Mid 1,n))\det(A) = \det(A(n\Mid n))\det(A(1\Mid 1)) - \det(A(n\Mid 1))\det(A(1\Mid n)). 
\end{equation*}
This can be written in the form
\begin{equation}
\det(A(1,n\Mid 1,n))\det(A) =
\det
\begin{pmatrix}
\det(A(1\Mid 1))\;&\;\det(A(1\Mid n))\\
\det(A(n\Mid 1))\;&\;\det(A(n\Mid n))
\end{pmatrix}.
\end{equation}
Numerous complicated identities involving determinants with determinants of submatrices as entries are developed in Muir~\cite{tm:ttd} Chapter VI, Compound Determinants.
\end{sremark}

%S. Gill Williamson, 2012\\
%${\rm cseweb.ucsd.edu\slash\sim gill\slash}$  
%
%\tableofcontents
%\index{contents}
%\hyperlink{index}{Index}
%\newpage
%\thispagestyle{empty}
\bibliographystyle{alpha}
\bibliography{AnchorLaplace}
\label{cor:relpridiaent}
\newpage
\mbox{}
\newpage
\hypertarget{index}{ }
\printindex
%\newpage
%\centerline{NOTES}
\end{document}